# UNIPOTENT REPRESENTATIONS: CHANGING $q$ TO $-q$

P. Deligne and G. Lusztig

ABSTRACT. Consider a Chevalley group over a finite field $F_q$ such that the longest element of the Weyl group is central. We construct an involution $\xi \mapsto \xi^!$ of the set of unipotent representations of this group such that the degree polynomial of a unipotent representation $\xi$ is obtained up to sign from the degree polynomial of $\xi^!$ by changing $q$ to $-q$.

## INTRODUCTION

**0.1.** Let **k** be an algebraic closure of a finite field $F_q$ with $q$ elements. Let $G$ be a reductive connected algebraic group over **k** with a fixed $F_q$-split rational structure with Frobenius map $F : G \to G$; let $G(F_q) = \{g \in G; F(g) = g\}$ (a finite group).

Let $W$ be the Weyl group of $G$, let $\lambda : W \to \mathbf{N}$ be the length function and let $w_0$ be the element of maximal length in $W$. In this paper we shall assume that $W$ is irreducible and that

(a) $w_0$ acts as $-1$ in the reflection representation of $W$.

For any finite group $\mathcal{G}$ let $\mathrm{Irr}(\mathcal{G})$ be the set of (isomorphism classes of) irreducible representations of $\mathcal{G}$ over $\bar{\mathbf{Q}}_l$, an algebraic closure of the field of $l$-adic numbers; here $l$ is a fixed prime number not dividing $q$.

Let $\mathcal{U}$ be the subset of $\mathrm{Irr}(G(F_q))$ consisting of unipotent representations (see [DL76, 7.8]). In [L84] a decomposition of $\mathrm{Irr}(W)$ into families is defined as well as a corresponding decomposition $\mathcal{U} = \sqcup_c \mathcal{U}_c$ indexed by the set $ce(W)$ of these families.

From [L84] it is known that $\mathcal{U}$ can be indexed in a way independent of $q$ (in the $E_7$ and $E_8$ cases a square root of $q$ and a square root of $-q$ are needed) and that for each $\xi \in \mathcal{U}$ the dimension $\dim(\xi)$ can be regarded as the value at $u = q$ of a polynomial $D_\xi(u) \in \mathbf{Q}[u]$ ($u$ is an indeterminate).

From the tables in [L84] one can see that for any $\xi \in \mathcal{U}_c$ there exists $\xi'$ in the same $\mathcal{U}_c$ such that $D_{\xi'}(u) = \pm D_\xi(-u)$. (This property has been used in [L25].) In general, $\xi'$ is not uniquely defined by this property. This happens for example for the unipotent representation $\xi$ in type $B_2$ with $D_\xi(u) = u(u^2+1)/2$; in this case $\xi'$ can be either $\xi$ or the other unipotent representation with $D_?(u) = u(u^2+1)/2$.

For $\xi \in \mathcal{U}$ let $A_\xi \in \mathbf{N}$ be the degree of the polynomial $D_\xi(u)$; let $a_\xi$ be the largest power of $u$ that divides $D_\xi(u)$.

Typeset by $\mathcal{A}\mathcal{M}\mathcal{S}$-TeX

**0.2.** In §2 we shall define an explicit involution $\mathcal{U} \to \mathcal{U}$, $\xi \mapsto \xi^!$, such that $D_\xi(-u) = (-1)^{A_\xi} D_{\xi^!}(u)$ and such that for any $c$, $\xi \mapsto \xi^!$ restricts to an involution $\mathcal{U}_c \to \mathcal{U}_c$.

In type $E_7, E_8$, in the cases where a square root of $q$ and one of $-q$ are needed for the parametrizations, this involution depends only on a choice of square root of $-1$ in $\bar{\mathbf{Q}}_l$.

In [L84] to each $c$ is attached a finite group $\Gamma_c$, a set $M(\Gamma_c)$ of simple objects in a category $Vec_{G_c}$ and an indexing of $\mathcal{U}_c$ by $M(\Gamma_c)$. The involution $\xi \mapsto \xi^!$ will be described using the category $Vec_{G_c}$.

**0.3.** Let $\mathcal{B}$ be the variety of Borel subgroups of $G$. It is defined over $F_q$. Let $F : \mathcal{B} \to \mathcal{B}$ be the corresponding Frobenius map. For a point of $\mathcal{B}$ identified with a Borel subgroup $B$ one has $(FB)(\mathbf{k}) = F(B(\mathbf{k}))$. It is known that $\mathcal{B} \times \mathcal{B}$ is a union of $G$-orbits (for the simultaneous $G$-conjugation action) which are parametrized by $W$. For $w \in W$ let $\mathcal{O}_w$ be the orbit indexed by $w$ and let $X_w = \{B \in \mathcal{B}; (B, F(B)) \in \mathcal{O}_w\}$ (see [DL76]). Note that $G(F_q)$ acts by conjugation on $X_w$ and this induces an action of $G(F_q)$ on the $l$-adic cohomology with compact support $H^i(X_w)$ (we omit the subscript indicating "compact"). This is a finite dimensional vector space over $\bar{\mathbf{Q}}_l$. Note that $H^*(X_w) = \sum_i (-1)^i H^i(X_w)$ can be viewed as an element of $\mathcal{K}$, the free abelian group with basis $\mathcal{U}$, see [DL76, 7.8].

Let $\mathcal{K}_{\mathbf{Q}} = \mathcal{K} \otimes \mathbf{Q}$. If $f \in \mathcal{K}_{\mathbf{Q}}$ we can write $f = \sum_{\xi \in \mathcal{U}} (\xi : f) \xi$ where $(\xi : f) \in \mathbf{Q}$. In §3 we shall prove the following result.

**Theorem 0.4.** *Let $w \in W$, $\xi \in \mathcal{U}$. We have*

$$(\xi : H^*(X_w)) = (-1)^{A_\xi}(\xi^!, H^*(X_{ww_0})).$$

A refinement of this theorem which takes weights into account is proved in §4 (see Theorem 4.3).

**0.5.** In §5 we describe (following [L78]) a connection between $X_{w_0}$ and the Iwahori-Hecke algebra with parameter $-q$ and we state some conjectural properties of it.

**0.6.** If property 0.1(a) is not satisfied, one can still define $\xi^!$ for $\xi \in \mathcal{U}$ with a property similar to that in 0.2; but this time $\xi^!$ is a unipotent representation of $G(F_q)$ where $G$ has an $F_q$-structure such that the Frobenius map acts on $W$ as conjugation by $w_0$. (See [E63] for $G$ of type $A$ and [L84, p.362] for $G$ of type $E_6$.)

1. THE ELEMENT $m_c \in M(\Gamma_c)$

**1.1.** Let $\Gamma$ be a finite group. Let $Vec_\Gamma$ be the category whose objects are finite dimensional $\bar{\mathbf{Q}}_l$-vector bundles on $\Gamma$ which are $\Gamma$-equivariant for the conjugation action of $\Gamma$ on $\Gamma$. For $U \in Vec_\Gamma$ let $U_g$ be the fibre of $U$ at $g \in \Gamma$. For $U, U' \in Vec_\Gamma$ we define the convolution $U * U' \in Vec_\Gamma$ by $(U * U')_g = \oplus_{g=ab} U_a \otimes U'_b$. This turns $Vec_\Gamma$ into a braided tensor category, see [L04]. Hence the Grothendieck group $K_\Gamma(\Gamma)$ of $Vec_G$ is a commutative ring.

For $U \in Vec_\Gamma$ and $g \in \Gamma$ the centralizer $Z(g)$ acts on $U_g$. For $U, U'$ in $Vec_\Gamma$ one defines
$$<U, U'> = |\Gamma|^{-1} \sum_{g,g'} \text{tr}(g, U'_{g'})\text{tr}(g'^{-1}, U_g),$$
the sum being extended to all commuting pairs $g, g'$ in $\Gamma$.

Let $M(\Gamma)$ be the set isomorphism classes of irreducible objects of $Vec_\Gamma$. Such an irreducible object $U$ has as support a single conjugacy class $C$ and, for $g \in C$, is determined by the representation $U_g$ of $Z(g)$, which should be irreducible. This identifies $M(\Gamma)$ with the set of $\Gamma$-conjugacy classes of pairs $(g, \rho)$ where $g \in \Gamma$, $\rho \in \text{Irr}(Z(g))$. We will simply write $(g, \rho)$ for the element of $M(\Gamma)$ defined by $(g, \rho)$. For $(g, \rho), (g', \rho')$ in $M(\Gamma)$ one has

$$<(g, \rho), (g', \rho')>$$
$$= |Z(g)|^{-1}|Z(g')|^{-1} \sum_{h \in \Gamma, h^{-1}ghg' = g'h^{-1}gh} \text{tr}(hg'h^{-1}, \rho)\text{tr}(h^{-1}g^{-1}h, \rho') \in \bar{\mathbf{Q}}_l,$$

see [L84, 4.14].

By [L87, 2.5], for any $(g', \rho') \in M_\Gamma$,

(a) *the group homomorphism $K_\Gamma(\Gamma) \to \bar{\mathbf{Q}}_l$,*
$(g, \rho) \mapsto (\dim \rho')^{-1}|Z(g')| <((g, \rho), (g', \rho')>$
*is a ring homomorphism.*

From (a) we see that for $x, y, z$ in $M(\Gamma)$ with $z = (g, \rho)$ we have

(b) $<z, x * y> = (\dim \rho)^{-1}|Z(g)| <z, x><z, y>$.

The invertible objects of the tensor category $Vec_\Gamma$ are the $(g, \rho)$ with $g$ central and $\rho$ a 1-dimensional representation of $\Gamma$. The set $M(\Gamma)_*$ of their isomorphism classes is an abelian group acting on $M(\Gamma)$ by convolution.

If $\Gamma$ is commutative then the abelian group $M(\Gamma)_* = M(\Gamma)$ is the abelian group $\Gamma \times \check{\Gamma}$ where $\check{\Gamma}$ is the Pontryagin dual of $\Gamma$. It is selfdual for a natural alternating selfduality and $<,>$ is the kernel of the corresponding Fourier transform normalized to be involutive.

If $\Gamma$ is the symmetric group $S_n, n \geq 3$ then $M(\Gamma)_*$ consists of $(1, \delta)$ where $\delta = 1$ or $\delta = \text{sgn}$ (the sign representation of $\Gamma$) and convolution with $(1, \delta)$ maps $(g, \rho)$ to $(g, \rho \otimes \delta|_{Z(g)})$.

**1.2.** For $E \in \text{Irr}(W)$ let $b_E$ be the smallest integer $\geq 0$ such that $E$ appears in the $b_E$-th symmetric power of the reflection representation of $W$.

The set $\text{Irr}(W)$ is partitioned into families (see [L84, §4]). Let $ce(W)$ be the set of families. We have $\text{Irr}(W) = \sqcup_{c \in ce(W)} c$. Let $c \in ce(W)$; let $\Gamma_c$ be the finite group attached to $c$ in [L84, §4].

Recall that we have an injective map $c \to M(\Gamma_c), E \mapsto m_E$, see [L84, §4]. Let $E(c) \in \text{Irr}(W)$ be the unique special representation in $c$. We have $m_{E(c)} = (1, 1)$. For $E \in c$ we set $b'_E = b_E - b_{E(c)}$. We have $b'_E \geq 0$.

**Theorem 1.3.** *If $|c| \neq 2$ (resp. $|c| = 2$), there is a unique (resp. exactly two) element(s) $m(c) \in M(\Gamma_c)_*$ such that for any $E \in c$ with $m_E = (g, \rho)$ we have*

(a) $<m_E, m(c)> = (-1)^{b'_E} \dim \rho / |Z(g)|$.

From the definitions in [L84,§4], one checks that if $|c| = 1$ so that $\Gamma_c = \{1\}$, there is nothing to prove (we have $1 = (-1)^0$). One can also check that, if $|c| = 2$ then $\Gamma_c = \mathbf{Z}/2$ so that $M(\Gamma_c) = \mathbf{Z}/2 \times \mathbf{Z}/2$. This case occurs only in type $E_7$ (one $c$, dimensions $512, 512$) or type $E_8$ (two $c$ exchanged by tensoring with sign, each with dimensions $4096, 4096$).

We now assume that $W$ is of type $B_n$ or $C_n$. As in [L84, §4], to $c$ we associate the symbol (up to shift) of the special representation $E(c)$. It is of the form

$$\Lambda = \begin{pmatrix} Z^2 \sqcup A \\ Z^2 \sqcup B \end{pmatrix}$$

where $Z^2, A, B$ are disjoint subsets of $\mathbf{N}$ with
$A = \{a_1 < a_2 < \cdots < a_{d+1}\}$, $B = \{b_1 < b_2 < \cdots < b_d\}$,
$a_1 < b_1 < a_2 < b_2 < \cdots < b_d < a_{d+1}$.
The sum of entries of this symbol is $n + (|Z^2| + d)^2$.

Let $Z^1 = A \cup B$. We have $|Z^1| = 2d+1$. Let $V$ be the $F_2$-vector space consisting of all subsets of $Z^1$ which have even cardinal with the sum of $U, U'$ being $U \sharp U'$. (For two subsets $U, U'$ of a set we define $U \sharp U' = (U \cup U') - (U \cap U')$.) We have $\dim V = 2d$. In this case, $\Gamma_c$ is abelian of order $2^d$ and $M(\Gamma_c)$ can be identified with the set of all symbols

$$\Lambda_Y = \begin{pmatrix} Z^2 \sqcup (Z^1 - Y) \\ Z^2 \sqcup Y \end{pmatrix}$$

for various subsets $Y \subset Z^1$ with $|Y| = d \mod 2$. Now $M(\Gamma_c)$ is in bijection with $V$ by $\Lambda_Y \leftrightarrow Y \sharp B$ and the sum operation on $M(\Gamma_c)$ corresponds to the sum operation on $V$. (See [L84, §4].) The pairing $<,>$ on $M(\Gamma_c)$ is given by

$$<\Lambda_Y, \Lambda_{Y'}> = 2^{-d}(-1)^{|(Y \sharp B) \cap (Y' \sharp B)|}.$$

The injective map $c \to M(\Gamma_c)$ has image consisting of all $\Lambda_Y$ such that $|Y| = d$; let $E_Y \in c$ be the representation corresponding to such a $\Lambda_Y$. We have $E(c) = E_B$. By [L79] we have

$$b'_{E_Y} = \sum_{b \in B} b - \sum_{y \in Y} y.$$

Hence

$$(-1)^{b'_{E_Y}} = (-1)^{|B \cap Z^1_{odd}| - |Y \cap Z^1_{odd}|}$$

where $Z^1_{odd} = Z^1 \cap \{1, 3, 5, \dots\}$. Let $Z^1_{ev} = Z^1 \cap \{0, 2, 4, \dots\}$. We have

$$|B \cap Z^1_{odd}| - |Y \cap Z^1_{odd}| + |B \cap Z^1_{ev}| - |Y \cap Z^1_{ev}| = |B \cap Z^1| - |Y \cap Z^1| = d - d = 0.$$

We set $Z^1_* = Z^1_{odd}$ if $|Z^1_{odd}| = 0 \mod 2$, $Z^1_* = Z^1_e v$ if $|Z^1_{ev}| = 0 \mod 2$. Note that $Z^1_*$ is well defined since $|Z^1|$ is odd. We see that

$$(-1)^{b'_{E_Y}} = (-1)^{|B \cap Z^1_*| - |Y \cap Z^1_*|} = (-1)^{|(Y \sharp B) \cap Z^1_*|}.$$

Thus

$$(-1)^{b'_{E_Y}} = 2^d < \Lambda_Y, \Lambda_{B \sharp Z^1_*} >.$$

We see that $m(c) = \Lambda_{B \sharp (Z^1_*)}$ satisfies (a). (We have $\dim \rho / |Z(g)| = 2^{-d}$ for any $(g, \rho) \in M(\Gamma_c)$.) To prove uniqueness it is enough to show that if $\Lambda_{Y'} \in M(\Gamma_c)$, $\Lambda_{Y''} \in M(\Gamma_c)$ satisfy $< \Lambda_Y, \Lambda_{Y'} > = < \Lambda_Y, \Lambda_{Y''} >$ for all $Y \subset Z^1$ such that $|Y| = d$ then $Y' = Y''$. This follows from the fact the subset $\{Y \sharp B; Y \subset Z^1, |Y| = d\}$ of $V$ contains a basis of $V$. This completes the proof in our case.

Next we assume that $W$ is of type $D_n$, $n$ even. We can assume that $|c| > 1$. As in [L84, §4], to $c$ we associate the symbol (up to shift) of the special representation $E(c)$. It is of the form

$$\Lambda = \begin{pmatrix} Z^2 \sqcup A \\ Z^2 \sqcup B \end{pmatrix}$$

where $Z^2, A, B$ are disjoint subsets of $\mathbf{N}$ with
$A = \{a_1 < a_2 < \cdots < a_d\}$, $B = \{b_1 < b_2 < \cdots < b_d\}$,
$b_1 < a_1 < b_2 < \cdots < b_d < a_d$.

The sum of entries of this symbol is $n + ((|Z^2| + d)^2 - (|Z^2| + d))$.
We have $d > 0$. Let $Z^1 = A \cup B$. We have $|Z^1| = 2d$ and $\sum_{z \in Z^1} z$ is even (since $n$ is even).

Let $V^+$ be the $F_2$-vector space consisting of all subsets of $Z^1$ which have even cardinal with the sum of $U, U'$ being $U \sharp U'$, taken modulo the line $\{0, Z^1\}$. We have $\dim V^+ = 2d - 2$. In this case, $\Gamma_c$ is abelian of order $2^{d-1}$ and $M(\Gamma_c)$ can be identified with the set of all symbols

$$\Lambda_Y = \Lambda_{\tilde{Y}} = \begin{pmatrix} Z^2 \sqcup \tilde{Y} \\ Z^2 \sqcup Y \end{pmatrix}$$

for various unordered pairs of subsets $Y \subset Z^1$, $\tilde{Y} \subset Z_1$ with $Z^1 = Y \sqcup \tilde{Y}$, $|Y| = d$ mod 2 (hence $|\tilde{Y}| = d \mod 2$). See [L84,§4]. Now $M(\Gamma_c)$ is in bijection with $V^+$ by $\Lambda_Y = \Lambda_{\tilde{Y}} \leftrightarrow Y \sharp B = \tilde{Y} \sharp B$ (equality modulo $Z^1$) and the sum operation on $M(\Gamma_c)$ corresponds to the sum operation on $V^+$. The pairing $<,>$ on $M(\Gamma_c)$ is given by

$$< \Lambda_Y, \Lambda_{Y'} > = 2^{-d+1} (-1)^{|(Y \sharp B) \cap (Y' \sharp B)|}.$$

The injective map $c \to M(\Gamma_c)$ has image consisting of all $\Lambda_Y = \Lambda_{\tilde{Y}}$ such that $|Y| = d$ (hence $|\tilde{Y}| = d$); let $E_Y = E_{\tilde{Y}} \in c$ be the representation corresponding to such a $\Lambda_Y = \Lambda_{\tilde{Y}}$. We have $E(c) = E_B$.

By [L79] we have

$b'_{E_Y} = \sum_{b \in B} b - \min(\sum_{y \in Y} y, \sum_{\tilde{y} \in \tilde{Y}} \tilde{y})$.

Hence $(-1)^{b'_{E_Y}}$ is either

$(-1)^{|B \cap Z^1_{odd}|}(-1)^{|Y \cap Z^1_{odd}|}$ or $(-1)^{|B \cap Z^1_{odd}|}(-1)^{|\tilde{Y} \cap Z^1_{odd}|}$

where $Z^1_{odd} = Z^1 \cap \{1, 3, 5, \dots\}$. We have

$\sum_{y \in Y} y + \sum_{\tilde{y} \in \tilde{Y}} \tilde{y} = \sum_{z \in Z^1} z$

and this is even. It follows that

$(-1)^{|Y \cap Z^1_{odd}|} = (-1)^{|\tilde{Y} \cap Z^1_{odd}|}$

hence

$(-1)^{b'_{E_Y}} = (-1)^{|B \cap Z^1_{odd}| + |Y \cap Z^1_{odd}|} = (-1)^{|B \cap Z^1_{odd}| + |\tilde{Y} \cap Z^1_{odd}|}$.

It follows that

$(-1)^{b'_{E_Y}} = (-1)^{|(B \sharp Y) \cap Z^1_{odd}|} = (-1)^{|(B \sharp \tilde{Y}) \cap Z^1_{odd}|}$.

Since $\sum_{z \in Z^1}$ is even we see that $|Z^1_{odd}|$ is even. Thus

$(-1)^{b'_{E_Y}} = 2^{d-1} < \Lambda_Y, \Lambda_{B \sharp Z^1_{odd}} >$.

We see that $m(c) = \Lambda_{B \sharp Z^1_{odd}}$ satisfies (a). (We have $\dim \rho / |Z(g)| = 2^{-d+1}$ for any $(g, \rho) \in M(\Gamma_c)$.) The proof of uniqueness is exactly like for type $B, C$. This completes the proof in our case.

In the $B, C, D$ cases, identifying $M(\Gamma_c)$ with a set of symbols, one checks that $m \mapsto m(c) * m$ corresponds to moving all odd entries in one row to the other row.

In the remainder of the proof we denote by $c_{odd}$ (resp. $c_{ev}$) the set of all $E \in c$ such that $b'_E$ is odd (resp. even); these $E$ can be determined from the tables in [BL78].

We now assume that $c$ is such that $\Gamma_c = S_5$. (In this case $W$ is of type $E_8$.) Then $c_{odd}$ (resp. $c_{ev}$) has 5 (resp. 12 elements) and the 5 (resp. 12) corresponding pairs $(g, \rho) \in M(S_5)$ are such that $g$ is an odd (resp. even) permutation. We see that for any $E \in c$ (with corresponding $(g, \rho) \in M(S_5)$) we have $(-1)^{b'_E} = \delta(g)$ where $\delta$ is the sign character of $S_5$. We have $\delta(g) = < (g, \rho), (1, \delta) > |Z(g)|/\dim \rho$ hence $m_c = (1, \delta)$ satisfies (a). The uniqueness is clear since $(1, 1)$ does not satisfy (a).

We now assume that $c$ is such that $\Gamma_c = S_4$. (In this case $W$ is of type $F_4$.) Then $c_{odd}$ (resp. $c_{ev}$) has 3 (resp. 8 elements) and the 3 (resp. 8) corresponding pairs $(g, \rho) \in M(S_4)$ are such that $g$ is an odd (resp. even) permutation. We see that for any $E \in c$ (with corresponding $(g, \rho) \in M(S_4)$) we have $(-1)^{b'_E} = \delta(g)$ where $\delta$ is the sign character of $S_4$. We have $\delta(g) = < (g, \rho), (1, \delta) > |Z(g)|/\dim \rho$ hence $m_c = (1, \delta)$ satisfies (a). The uniqueness is clear since $(1, 1)$ does not satisfy (a).

We now assume that $c$ is such that $\Gamma_c = S_3$. (In this case $W$ is of type $E_7, E_8$ or $G_2$.) If $c_{odd}$ is empty (which can happen in type $E_8$) then $c = c_{ev}$ has 5 elements. In this case for any $E \in c$ (with corresponding $(g, \rho) \in M(S_3)$) we have $(-1)^{b'_E} = 1 = < (g, \rho), (1, 1) > |Z(g)|/\dim \rho$ hence $m_c = (1, 1)$ satisfies (a). The uniqueness is clear since $(1, \text{sgn})$ does not satisfy (a).

We now assume that $c_{odd}$ is nonempty. Then $c_{odd}$ has 1 element and $c_{ev}$ has 4 elements (in type $E_7, E_8$) and 3 elements (in type $G_2$). The element $(g, \rho)$ corresponding to $E \in c_{odd}$ is such that $g$ is an odd permutation; the elements $(g, \rho)$ corresponding to the various $E \in c_{ev}$ are such that $g$ is an even permutation. We see that for any $E \in c$ (with corresponding $(g, \rho) \in M(S_3)$) we have $(-1)^{b'_E} = \delta(g)$ where $\delta$ is the sign character of $S_3$. We have $\delta(g) = <(g, \rho), (1, \delta)> |Z(g)|/\dim \rho$ hence $m_c = (1, \delta)$ satisfies (a). The uniqueness is clear since $(1, 1)$ does not satisfy (a).

We now assume that $c$ is such that $\Gamma_c = S_2$, $|c| > 2$ and that $W$ is of type $E_6, E_7, E_8$ or $F_4$. In this case $M(\Gamma_c)$ is an $F_2$-vector space of dimension 2. We have $|c| = 3$ so that the image of $c$ in this vector space is a basis of the vector space union with 0. Hence there is a unique $F_2$-linear form on this vector space whose composition with $c \to M(\Gamma_c)$ takes the value $b'_E \mod 2$ at any $E \in c$. This linear form must be of the form $m \mapsto 2 < m, m(c) >$ for a well defined $m(c) \in M(\Gamma_c) = M(\Gamma_c)_*$. Thus the theorem holds in this case.

We now assume that $\Gamma_c = S_2, |c| = 2$. Again, $M(\Gamma_c)$ is an $F_2$-vector space of dimension 2 and the image of $c$ in this vector space is a subspace of dimension 1. Hence there are exactly two $F_2$-linear forms on this vector space whose composition with $c \to M(\Gamma_c)$ take the value $b'_E \mod 2$ at any $E \in c$. Such a linear form must be of the form $m \mapsto 2 < m, m(c) >$ for a well defined $m(c) \in M(\Gamma_c) = M(\Gamma_c)_*$. Thus the theorem holds in this case. This completes the proof.

**Proposition 1.4.** *Let $m(c)$ be as in 1.3. For any $E \in c$, $m \in M(\Gamma_c)$ we have*

$$< m_E, m(c) * m > = (-1)^{b'_E} < m_E, m > .$$

Using 1.1(b) we see that this follows from 1.3(a).

## 2. The involution $\xi \mapsto \xi^!$

**2.1.** For any $E \in \text{Irr}(W)$ we define
$R_E = |W|^{-1} \sum_{w \in W} \text{tr}(w, E) H^*(X_w) \in \mathcal{K}_{\mathbf{Q}}$.
For $c \in ce(W)$ let $\mathcal{U}_c$ be the set of all $\xi \in \mathcal{U}$ such that $(\xi : R_E) \neq 0$ for some $E \in c$. By [L84, 4.23] we have $\mathcal{U} = \sqcup_{c \in ce(W)} \mathcal{U}_c$ and for any $c \in ce(W)$ we have a bijection $M(\Gamma_c) \xrightarrow{\sim} \mathcal{U}_c$, $m \mapsto \xi_m$.

We have the following result.

**Theorem 2.2.** *Let $c \in ce(W)$. Let $m(c)$ be as in 1.3 and let $m \in M(\Gamma_c)$. The polynomials $D_{\xi_m}(u)$ and $D_{\xi_{m(c)*m}}(u)$ giving the dimensions of the unipotent representations $\xi_m$ and $\xi_{m(c)*m}$ corresponding to $m$ and $m(c) * m$, are deduced one from the other, up to sign, by replacing $u$ by $-u$.*

For $c \in ce(W)$ let $a_c = a_\xi$, $A_c = A_\xi$ where $\xi \in \mathcal{U}_c$. (This is well defined.)

Define $\Delta : M(\Gamma_c) \to \{1, -1\}$ to be identically 1 if $|c| \neq 2$ while if $|c| = 2$ we define $\Delta$ to be 1 on the image of $c \to M(\Gamma_c)$ and to be $-1$ on the complement of that image. From the definition we have

(a) $\Delta(m(c) * m)\Delta(m) = (-1)^{a_c + A_c}$.

For $E \in \text{Irr}(W)$ we can define $\dim(R_E)(q)$ to be

$|W|^{-1} \sum_{w \in W} \text{tr}(w, E) \sum_i (-1)^i \dim H^i(X_w)$.

By [L78, (3.19.1)] we have

$\dim(R_E)(q) = \sum_{i \in \mathbf{N}} (E : \bar{S}_i) q^i$

where $(E : \bar{S}_i)$ is the multiplicity of $E$ in a certain $W$-module $\bar{S}_i$ (which is 0 for large $i$); in particular, $\dim(R_E)(q)$ is a polynomial in $q$ with coefficients in $\mathbf{N}$. We show that

$\dim(R_E)(-q) = (-1)^{b_E} \dim(R_E(q))$.

It is enough to show that $(E : \bar{S}_i) = 0$ for any $i$ such that $i \neq b_E \mod 2$; indeed, using 0.1(a) and the definition of $\bar{S}_i$, we see that $w_0$ acts on $\bar{S}_i$ as multiplication by $(-1)^i$ and on $E$ as multiplication by $(-1)^{b_E}$.

From [DL76, 7.9] and [L84, 4.23] we have

(b) $$D_{\xi_m}(q) = \sum_{E \in c} \Delta(m) <m_E, m> \dim(R_E)(q).$$

Replacing $q$ by $-q$ we obtain:

(c) $$D_{\xi_m}(-q) = \sum_{E \in c} \Delta(m) <m_E, m> \dim(R_E)(-q).$$

Replacing in (b) $m$ by $m(c) * m$ we obtain

(d) $$D_{\xi_{m(c)*m}}(q) = \sum_{E \in c} \Delta(m(c) * m) <m_E, m(c) * m> \dim(R_E)(q).$$

Using 1.4 we obtain

$$D_{\xi_{m(c)*m}}(q) = \sum_{E \in c} \Delta(m(c) * m)(-1)^{b_{E(c)}}(-1)^{b_E} <m_E, m> \dim(R_E)(q)$$

$$= \Delta(m(c) * m)(-1)^{b_{E(c)}} \sum_{E \in c} <m_E, m> \dim(R_E)(-q).$$

Using (c) we obtain

$$D_{\xi_{m(c)*m}}(q) = \Delta(m(c) * m)\Delta(m)(-1)^{b_{E(c)}} D_{\xi_m}(-q)$$

that is (using (a) and $b_{E(c)} = a_c$):

$$D_{\xi_{m(c)*m}}(q) = (-1)^{A_c} D_{\xi_m}(-q).$$

This proves the theorem.

**2.3.** Let $c \in ce(W)$ be such that $|c| \neq 2$. Then $m(c)$ in 1.3 is uniquely determined. We define a bijection $M(\Gamma_c) \to M(\Gamma_c)$, $m \mapsto m^!$, by $m^! = m(c) * m$. We define $\mathcal{U}_c \to \mathcal{U}_c, \xi \mapsto \xi^!$ by $(\xi_m)^! = \xi_{m^!}$ for any $m \in M(\Gamma)$. Since $m' * m' = (1,1)$ for any $m' \in M(\Gamma)_*$, we see that $\xi \mapsto \xi^!$ is an involution of $\mathcal{U}_c$.

**2.4.** Let $c \in ce(W)$ be such that $|c| = 2$. (In this case $W$ is of type $E_7$ or $E_8$). We have $\mathcal{U}_c = \mathcal{U}'_c \sqcup \mathcal{U}''_c$ where $\mathcal{U}'_c = \{\xi_{m_E}; E \in c\}$ and $\mathcal{U}''_c = \mathcal{U}_c - \mathcal{U}'_c$. While the map $c \to \mathcal{U}'_c$, $E \mapsto \xi_{m_E}$ depends on a choice of $\sqrt{q}$, its image $\mathcal{U}'_c$ does not. Now $\mathcal{U}''_c$ consists of two unipotent representations whose associated eigenvalue of Frobenius (see [L78, 3.9] and [L84, §11]) is a square root of $-q$ times an integer power of $q$, for one of them, and minus that square of $-q$ times an integer power of $q$, for the other one.

Now there are two choices of $m(c)$ in 1.3 hence two bijections $\mathcal{U}_c \to \mathcal{U}_c$, $\xi_m \mapsto \xi_{m(c)+m}$ which restrict to the two possible bijections $\mathcal{U}'_c \to \mathcal{U}''_c$.

Assume now that we have fixed $\sqrt{-1} \in \bar{\mathbf{Q}}_l$. We define a bijection $\mathcal{U}'_c \to \mathcal{U}''_c$ by the requirement that it takes $\xi_{(1,1)}$ (defined using a choice of $\sqrt{q}$) to the unipotent representation in $\mathcal{U}''_c$ with associated eigenvalue of Frobenius $\sqrt{-1}\sqrt{q}$ times an integer power of $q$. This bijection $\mathcal{U}'_c \to \mathcal{U}''_c$ is independent of the choice of $\sqrt{q}$ and $\sqrt{-q}$ (but it depends on the choice of $\sqrt{-1}$). We define $\mathcal{U}_c \to \mathcal{U}_c, \xi \mapsto \xi^!$ to be the bijection which restricts to $\mathcal{U}'_c \to \mathcal{U}''_c$ as above. This bijection is of the form $\xi_m \to \xi_{m(c)*m}$ (with $m \in \Gamma_c$) for one particular choice of $m(c)$ (determined by $\sqrt{-1}$). We also define a bijection $M(\Gamma_c) \to M(\Gamma_c)$, $m \mapsto m^!$, by $m^! = m(c) * m$.

**2.5.** We define a bijection $\mathcal{U} \to \mathcal{U}$, $\xi \mapsto \xi^!$ which for any $c \in ce(W)$ with $|c| \neq 2$ restricts to the bijection $\mathcal{U}_c \to \mathcal{U}_c$ in 2.3 and for any $c \in ce(W)$ with $|c| = 2$ restricts to the bijection $\mathcal{U}_c \to \mathcal{U}_c$ in 2.4. This bijection is as in 0.2. It is an involution.

## 3. Proof of Theorem 0.4

**3.1.** Let $c \in ce(W)$, $m \in M(\Gamma_c)$. We have $H^*(X_w) = \sum_{E \in \mathrm{Irr}(W)} \mathrm{tr}(w, E) R_E$. Hence

$$(\xi_m : H^*(X_w)) = \sum_{E \in \mathrm{Irr}(W)} \mathrm{tr}(w, E)(\xi_m : R_E) = \sum_{E \in c} \mathrm{tr}(w, E)(\xi_m : R_E).$$

By [L84, 4.23] for any $E \in c$ we have

(a) $$(\xi_m : R_E) = \Delta(m) <m_E, m>.$$

Hence

$$(\xi_m : H^*(X_w)) = \sum_{E \in c} \mathrm{tr}(w, E) \Delta(m) <m_E, m>.$$

It is then enough to show that

$$\sum_{E \in c} \mathrm{tr}(w, E) \Delta(m) <m_E, m> = (-1)^{A_c} \sum_{E \in c} \mathrm{tr}(ww_0, E) \Delta(m^!) <m_E, m^!>$$

or that

$$\operatorname{tr}(w, E)\Delta(m) <m_E, m> = (-1)^{A_c}\operatorname{tr}(ww_0, E)\Delta(m^!) <m_E, m^!>$$

for any $E \in c$. Thus we must show that

$$\operatorname{tr}(w, E)\Delta(m) <m_E, m> = (-1)^{A_c}\operatorname{tr}(ww_0, E)\Delta(m^!)(-1)^{b'_E} <m_E, m>$$

for any $E \in c$. Using 0.1(a) and the definition of $b_E$ we see that
  (b) $w_0$ *acts on $E$ as multiplication by $(-1)^{b_E}$.*
Thus

$$\operatorname{tr}(ww_0, E) = (-1)^{b_E}\operatorname{tr}(w, E)$$

and it is enough to show that $\Delta(m) = \Delta(m^!)(-1)^{b_{E(c)}+A_c}$. This follows from 2.2(a) using $b_{E(c)} = a_c$. This proves Theorem 0.4.

**3.2.** Let $G(F_q)_{rs}$ be the set of regular semisimple elements in $G(F_q)$. Let $cl(W)$ be the set of conjugacy classes in $W$. For $s \in G(F_q)_{rs}$ we choose $B \in \mathcal{B}$ such that $s \in B$. Then $B \in X_w$ where the conjugacy class $(w)$ of $w \in W$ is independent of the choice of $B$. Thus $s \mapsto (w)$ is a well defined map $\tau : G(F_q)_{rs} \to cl(W)$.

We now assume that $q$ is large enough so that $\tau$ is surjective.

Let $\gamma \in cl(W)$; note that $w_0\gamma \in cl(W)$. Let $s \in \tau^{-1}(\gamma)$, $s' \in \tau^{-1}(w_0\gamma)$. We show:
  (a) *For any $\xi \in \mathcal{U}$ we have $\operatorname{tr}(s, \xi) = (-1)^{A_\xi}\operatorname{tr}(s', \xi^!)$.*
By [DL76, 7.9] we have $\operatorname{tr}(s, \xi) = (\xi, H^*(X_w))$ and similarly $\operatorname{tr}(s', \xi^!) = (\xi^!, H^*(X_{ww_0}))$ where $w \in \gamma$. Hence (a) follows from Theorem 0.4.

## 4. Weight spaces

**4.1.** Let $v$ be an indeterminate. Let $\mathbf{H}$ be the Iwahori-Hecke algebra attached to $W$. This is a free $\bar{\mathbf{Q}}_l[v]$-module with basis $\{T_w; w \in W\}$; the product is defined by $T_w T_{w'} = T_{ww'}$ if $\lambda(ww') = \lambda(w) + \lambda(w')$ and $(T_s + 1)(T_s - v^2) = 0$ if $s \in W$, $\lambda(s) = 1$. Let $\underline{\mathbf{H}} = \bar{\mathbf{Q}}_l(v) \otimes_{\bar{\mathbf{Q}}_l[v]} \mathbf{H}$; this is an algebra over $\bar{\mathbf{Q}}_l(v)$ canonically isomorphic to the group algebra of $W$ over $\bar{\mathbf{Q}}_l(v)$, see [L81]. For $E \in \operatorname{Irr}(W)$ we denote by $E(v)$ the corresponding simple $\underline{\mathbf{H}}$-module. From 0.1(a) we see that $T_{w_0}$ is in the centre of $\underline{\mathbf{H}}$; it follows that if $E \in c$, $T_{w_0}$ acts on $E(v)$ as multiplication by an element in $\bar{\mathbf{Q}}_l(v)$ which, by [L84,(5.12.2)], is of the form $ev^{2\nu-a_c-A_c}$ where $\nu = \lambda(w_0)$ and $e = \pm 1$ is such that the action of $w_0$ on $E$ is by multiplication by $e$. By 3.1(b) we have $e = (-1)^{b_E}$.

**4.2.** Let $w \in W$. Using [D80] we have $H^i(X_w) = \oplus_{k \in \mathbf{Z}} H^i_k(X_w)$ where $H^i_k(X_w)$ is the part of weight $k$ of $H^i(X_w)$. (this is a $G^F$-submodule).

By [DM78], [A83] (see also [L84, 3.8]), for $w \in W$ we have

(a) $$\sum_{E \in \operatorname{Irr}(W)} \operatorname{tr}(T_w, E(v)) R_E = \sum_{k \in \mathbf{Z}} H^*_k(X_w) v^k$$

where $H_k^*(X_w) = \sum_i (-1)^i H_k^i(X_w)$; in particular we have $H_k^*(X_w) \in \mathcal{K}$; (a) is an equality in $\mathcal{K} \otimes \mathbf{Q}[v]$. Let $c \in ce(W)$, $m \in M(\Gamma_c)$. From (a) we have

(b) $$\sum_{k \in \mathbf{Z}} (\xi_m : H_k^*(X_w)) v^k = \sum_{E \in c} \mathrm{tr}(T_w, E(v))(\xi_m : R_E).$$

Using 3.1(a) we deduce

$$\sum_{k \in \mathbf{Z}} (\xi_m : H_k^*(X_w)) v^k = \sum_{E \in c} \mathrm{tr}(T_w, E(v)) \Delta(m) <m_E, m>.$$

Replacing $w, m$ by $ww_0, m^!$ and using

$$\mathrm{tr}(T_{ww_0}, E(v)) = \mathrm{tr}(T_{w^{-1}} T_{w_0}, E(v)) = (-1)^{b_E} v^{2\nu - a_c - A_c} \mathrm{tr}(T_{w^{-1}}, E(v))$$
$$= (-1)^{b_E} v^{2\nu - a_c - A_c} \mathrm{tr}(T_w, E(v))$$

(which follows from 4.1) we obtain

$$\sum_{k \in \mathbf{Z}} (\xi_{m^!} : H_k^*(X_{ww_0})) v^k = \sum_{E \in c} (-1)^{b_E} v^{2\nu - a_c - A_c} \mathrm{tr}(T_w, E(v)) \Delta(m^!) <m_E, m^!>.$$

Using $<m_E, m^!> = (-1)^{b'_E} <m_E, m^!>$ (see 1.4) and $\Delta(m^!) = (-1)^{a_c + A_c} \Delta(m)$ see 2.2(a) we deduce

$$\sum_{k \in \mathbf{Z}} (\xi_{m^!} : H_k^*(X_{ww_0})) v^k = (-1)^{A_c} \sum_{E \in c} v^{2\nu - a_c - A_c} \mathrm{tr}(T_w, E(v)) \Delta(m) <m_E, m>.$$

Using now (b) we deduce

$$\sum_{k \in \mathbf{Z}} (\xi_{m^!} : H_k^*(X_{ww_0})) v^k = (-1)^{A_c} v^{2\nu - a_c - A_c} \sum_{k \in \mathbf{Z}} (\xi_m : H_k^*(X_w)) v^k$$
$$= (-1)^{A_c} \sum_{k \in \mathbf{Z}} (\xi_m : H_{k - 2\nu + a_c + A_c}^*(X_w)) v^k.$$

Hence for any $k \in \mathbf{Z}$ the following refinement of 0.4 holds.

**Theorem 4.3.**
$$(\xi_{m^!} : H_k^*(X_{ww_0})) = (-1)^{A_c} (\xi_m : H_{k - 2\nu + a_c + A_c}^*(X_w)).$$

It follows that

(a) $$H_k^*(X_{ww_0}) = \sum_{\xi \in \mathcal{U}} (-1)^{A_\xi} (\xi : H_{k - 2\nu + a_\xi + A_\xi}^*(X_w)) \xi^!.$$

Taking sum over $k$ we obtain

(b) $$H^*(X_{ww_0}) = \sum_{\xi \in \mathcal{U}} (-1)^{A_\xi} (\xi : H^*(X_w)) \xi^!.$$

## 5. The Iwahori-Hecke algebra with parameter $-q$

**5.1.** Let $\mathbf{H}_{\sqrt{q}}$ (resp. $\mathbf{H}_{\sqrt{-q}}$) be the specialization of $\mathbf{H}$ for $v = \sqrt{q}$ (resp. $v = \sqrt{-q}$). This is an algebra over $\bar{\mathbf{Q}}_l$. The method of [L81] gives an algebra isomorphism of $\mathbf{H}_{\sqrt{q}}$ (resp. $\mathbf{H}_{\sqrt{-q}}$) with the group algebra of $W$ over $\bar{\mathbf{Q}}_l$. For $E \in \mathrm{Irr}(W)$ we denote by $E(\sqrt{q})$ (resp. $E(\sqrt{-q})$) the corresponding simple $\mathbf{H}_{\sqrt{q}}$-module (resp. $\mathbf{H}_{\sqrt{-q}}$-module).

**5.2.** Now $H^*(X_1) = H^0(X_1)$ is naturally an $\mathbf{H}_{\sqrt{q}}$-module which commutes with the $G(F_q)$-action. We can decompose $H^*(X_1)$ with respect to the action of $\mathbf{H}_{\sqrt{q}}, G(F_q)$ as $\oplus_{E \in \mathrm{Irr}(W)} E(\sqrt{q}) \otimes \xi_E$ where $\xi_E \in \mathcal{U}$. Also $E \mapsto \xi_E$ is an imbedding of $\mathrm{Irr}(W)$ into $\mathcal{U}$. .

**5.3.** Following [L78, p.24], for any $w \in W$ we define a map $t_w : X_{w_0} \to X_{w_0}$ by $B \mapsto B'$ where $B' \in \mathcal{B}$ is defined by the conditions

$$(B, B') \in \mathcal{O}_w, (B', F(B)) \in \mathcal{O}_{w^{-1}w_0}.$$

In *loc.cit.* it is shown that the linear maps $(-1)^{\lambda(w)} t_w^*$ on $H^i(X_{w_0})$ define a representation of the Iwahori-Hecke algebra $\mathbf{H}_{\sqrt{-q}}$. We can decompose $H^i(X_{w_0})$ with respect to the action of $\mathbf{H}_{\sqrt{-q}}, G(F_q)$ as

$$\oplus_{E \in \mathrm{Irr}(W)} E(\sqrt{-q}) \otimes H^i(X_{w_0})_E$$

where $H^i(X_{w_0})_E$ is a $G(F_q)$-module. It follows that we have

$$H^*(X_{w_0}) = \sum_{E \in \mathrm{Irr}(W)} \dim(E) H^*(X_{w_0})_E$$

(equality of virtual representations of $G(F_q)$) where

$$H^*(X_{w_0})_E = \sum_i (-1)^i H^i(X_{w_0})_E.$$

**5.4.** From 4.3(b) with $w = 1$ we have

$$H^*(X_{w_0}) = \sum_{E \in \mathrm{Irr}(E)} (-1)^{A_{\xi_E}} \dim(E) \xi_E^!.$$

More precisely from 4.3(a) with $w = 1$ we have for any $k \in \mathbf{Z}$:

$$H_k^*(X_{w_0}) = \sum_{E \in \mathrm{Irr}(W); 2\nu - a_{\xi_E} - A_{\xi_E} = k} (-1)^{A_{\xi_E}} \dim(E) \xi_E^!.$$

(We use that $H^*(X_1) = H_0^*(X_1)$.)

We conjecture that for any $E$ in the family $c \subset \text{Irr}(W)$ we have

$$H^*(X_{w_0})_E = (-1)^{A_c}(\xi_E)^!.$$

In the case where $W$ is of type $E_7$ or $E_8$ so that ! is defined using a choice of $\sqrt{-1}$, we assume that $\sqrt{-q} = \sqrt{-1}\sqrt{q}$. The same assumption is made in the following stronger conjecture.

$$H^{2\nu-A_c}(X_{w_0})_E = (\xi_E)^!,$$

$$H^i(X_{w_0})_E = 0 \text{ if } i \neq 2\nu - A_c.$$

We also expect that $H^{2\nu-A_c}(X_{w_0})_E$ is contained in $H^{2\nu-A_c}_{2\nu-a_c-A_c}(X_{w_0})$.

Institute for Advanced Study, Princeton, NJ 08540; Department of Mathematics, M.I.T., Cambridge, MA 02139